\documentclass[pdflatex,sn-mathphys-num]{sn-jnl} 

\usepackage{graphicx}
\usepackage{multirow}
\usepackage{amsmath,amssymb,amsfonts}
\usepackage{amsthm}
\usepackage{mathrsfs}
\usepackage[title]{appendix}
\usepackage{xcolor}
\usepackage{textcomp}
\usepackage{manyfoot}
\usepackage{booktabs}
\usepackage{algorithm}
\usepackage{algorithmicx}
\usepackage{algpseudocode}
\usepackage{listings}

\theoremstyle{thmstyleone}
\newtheorem{theorem}{Theorem}

\theoremstyle{thmstyletwo}

\theoremstyle{thmstylethree}

\newcommand{\sgn}{\operatorname{sgn}}

\begin{document}

\title[On the CDF of a Correlated Gaussian Product (Zero Means)]{On the fully analytical cumulative distribution of product of correlated Gaussian random Variables with zero means}

\author*[1]{\fnm{Asst. Prof. Erdinc} \sur{Akyildirim, Phd}}\email{erdinc.akyildirim@nottingham.ac.uk}
\author[2]{\fnm{Alper} \sur{Hekimoglu, Phd}}\email{a.hekimoglu@eib.org}
\equalcont{These authors contributed equally to this work.}

\affil*[1]{\orgdiv{Department of Finance, Accounting and Banking}, \orgname{University of Nottingham}, \orgaddress{\street{Street}, \city{Nottingham}, \postcode{100190},  \country{United Kingdom}}}
\affil[2]{\orgdiv{Model Validation Unit}, \orgname{European Investment Bank}, \orgaddress{\street{98-100, boulevard Konrad Adenauer}, \city{Luxembourg City}, \postcode{L-2950}, \country{Luxembourg}}}

\abstract{
We derive a fully analytical, one-line closed-form expression for the cumulative distribution function (CDF) of the product of two correlated zero-mean normal random variables, avoiding any series representation. This result complements the well-known compact density formula with an equally compact and computationally practical CDF representation.

Our main formula expresses the CDF in terms of Humbert’s confluent hypergeometric function $\Phi_1$ and modified Bessel functions $K_\nu$, offering both theoretical elegance and computational efficiency. High-precision numerical experiments confirm pointwise agreement with Monte Carlo simulations and other benchmarks to machine accuracy.

The resulting representation provides a tractable tool for applications in wireless fading channel modeling, nonlinear signal processing, statistics, finance, and applied probability.

}


\maketitle

\section{Introduction}\label{sec1}

This is a short note on the cumulative distribution function of the product of bivariate normal random variables with correlation $\rho$. We also derive distributions the sum and mean of the these variables with sample size $n$. Our main purpose is to derive a representation that is fully analytical, meaning either a formula composed fully of elementary functions or special functions widely used in the literature.  

Closed-form expressions for the product of correlated normal variables have a long history. 
Craig \cite{Craig1936} first analyzed the independent case $\rho=0$ and obtained density 
representations involving modified Bessel functions. 
Nadarajah et al.\ \cite{NadarajahPogany2016} extended these results to the general correlated case with zero means, 
providing an explicit formula for the density in terms of the modified Bessel function $K_0$. 
Cui et al.\ \cite{Cui2016} further analyzed the non-central correlated case, deriving infinite-series 
representations for the density.

However, recent reviews \cite{Gaunt2021,Gaunt2022} explicitly state that 
a closed-form expression for the cumulative distribution function (CDF) of the product 
of two correlated normal variables was not previously available, except in special cases 
(e.g., $\rho=0$) where the CDF can be expressed in terms of Bessel or Struve functions. However, \cite{Gaunt2019ProductDistribution} proved the link between this product distibution and variance gamma distribution, which is an important link to have alternative representations. 

More recently, Gaunt \cite{Gaunt2023VG} derived an exact CDF representation for the Variance--Gamma distribution---which includes the product of zero-mean correlated normals as a special case---in terms of infinite series involving modified Bessel and Lommel functions. 
Though elegant, these series-form CDF formulas may require special software and truncation handling. 
Our result directly fills this gap, and can be viewed as the \emph{CDF analogue} of the closed-form 
density representation given by \cite{NadarajahPogany2016}. Our closed form leads to fast evaluation, proven machine-precision accuracy, and easier extension to sample means and sums, as we demeonstrate in our numerical experiments.

Therefore, to the best of our knowledge, the one-line closed-form CDF formula 
involving Humbert's Confluent Hypergeometric function and modified Bessel function of the second kind presented in this paper 
appears to be the first of its kind in the literature.

\section{Main Results}\label{sec4}

\begin{theorem}\label{thm:theorem1}
Let $(X,Y)$ denote a bivariate normal random vector with zero means, variances $\sigma_x,\sigma_y$ and correlation coefficient $\rho$. Let $\rho\in(-1,1)$ and $z\in\mathbb{R}$. Then the CDF can be written as
\begin{align}
F(z\mid \rho)
\label{eq:analytical_cdf_humbert}
&=\frac{C\,e^{sC}\,\sqrt{1-\rho}}{\sqrt{2}\,\pi}\Bigg[
\;2\,K_{1}(C)\,\Phi_{1}\!\Big(\tfrac12,\tfrac12;\tfrac32;\ \kappa,\ Y\Big)\nonumber\\
&\quad-\;s\,\tfrac{2}{3}(1-\rho)\,K_{0}(C)\,\Phi_{1}\!\Big(\tfrac32,\tfrac12;\tfrac52;\ \kappa,\ Y\Big)\nonumber\\
&\quad+\;2s\,K_{0}(C)\,\Phi_{1}\!\Big(\tfrac12,\tfrac12;\tfrac32;\ \kappa,\ Y\Big)
\Bigg],
\end{align}
where $K_{\nu}$ is the modified Bessel function of the second kind and
$\Phi_{1}$ is Humbert’s confluent hypergeometric function. Use $C>0$
in $K_{\nu}$ and in the power; the signed $sC$ appears only in $e^{sC}$
and in $Y$. By continuity, $F(0\mid\rho)=\tfrac12-\frac{\arcsin{\rho}}{\pi}$. Define
\[
B:=1-\rho^{2},\qquad C:=\frac{|z|}{B},\qquad s:=\sgn(z),\qquad
\kappa:=\frac{1-\rho}{2},\qquad Y:=-\,s\,C\,(1-\rho).
\]
\end{theorem}
\begin{theorem}\label{thm:theorem2}
Let $(Z_1,Z_2,\cdot\cdot\cdot Z_{\hat{n}})$ is distribued according to \ref{thm:theorem2}. Let $\rho\in(-1,1)$, $\hat{Z}$ denote their sample mean and  $z\in\mathbb{R}$. Then the CDF can be written as
\begin{align}
F(z\mid \rho)
\label{eq:analytical_cdf_humbertMean}
&=\frac{C\,e^{sC}\,\left(\sqrt{1-\rho}\right)^{\hat{n}}}{\sqrt{2\pi}\,\Gamma(\hat{n}/2)}\Bigg[
\;\frac{2}{n}\,K_{\tfrac{\hat{n}+1}{2}}(C)\,\Phi_{1}\!\Big(\tfrac {\hat{n}}{2},1-\tfrac {n}{2};\tfrac {\hat{n}}{2}+1 ;\ \kappa,\ Y\Big)\nonumber\\
&\quad-\;s\,\frac{2}{2+\hat{n}}(1-\rho)\,K_{\tfrac{\hat{n}-1}{2}}(C)\,\Phi_{1}\!\Big(\tfrac{ \hat{n}}{2}+1;1-\tfrac{\hat{n}}{2};\tfrac {\hat{n}}{2}+2;\ \kappa,\ Y\Big)\nonumber\\
&\quad+\;\frac{2}{\hat{n}}s\,K_{\tfrac{\hat{n}-1}{2}}(C)\,\Phi_{1}\!\Big(\tfrac {\hat{n}}{2},1-\tfrac {\hat{n}}{2};\tfrac {\hat{n}}{2}+1;\ \kappa,\ Y\Big)
\Bigg],
\end{align}
where $K_{\nu}$ is the modified Bessel function of the second kind and
$\Phi_{1}$ is Humbert’s confluent hypergeometric function. Use $C>0$
in $K_{\nu}$ and in the power; the signed $sC$ appears only in $e^{sC}$
and in $Y$. By continuity, $F(0\mid\rho)=\tfrac12-\frac{\arcsin{\rho}}{\pi}$. Define
\[
B:=1-\rho^{2},\qquad C:=\frac{\hat{n}|z|}{B},\qquad s:=\sgn(z),\qquad
\kappa:=\frac{1-\rho}{2},\qquad Y:=-\,s\,C\,(1-\rho).
\]
\end{theorem}

\begin{theorem}\label{thm:theorem3}
Let $(Z_1,Z_2,\cdot\cdot\cdot Z_{\hat{n}}$  is distribued according to \ref{thm:theorem2}. Let $\rho\in(-1,1)$, $Z_{\Sigma}$ denote their sample sum and  $z\in\mathbb{R}$. Then the CDF can be written as
\begin{align}
F(z\mid \rho)
\label{eq:analytical_cdf_humbertSum}
&=\frac{C\,e^{sC}\,\left(\sqrt{1-\rho}\right)^{\hat{n}}}{\sqrt{2\pi}\,\Gamma(\hat{n}/2)}\Bigg[
\;\frac{2}{n}\,K_{\tfrac{\hat{n}+1}{2}}(C)\,\Phi_{1}\!\Big(\tfrac {\hat{n}}{2},1-\tfrac {n}{2};\tfrac {\hat{n}}{2}+1 ;\ \kappa,\ Y\Big)\nonumber\\
&\quad-\;s\,\frac{2}{2+\hat{n}}(1-\rho)\,K_{\tfrac{\hat{n}-1}{2}}(C)\,\Phi_{1}\!\Big(\tfrac{ \hat{n}}{2}+1;1-\tfrac{\hat{n}}{2};\tfrac {\hat{n}}{2}+2;\ \kappa,\ Y\Big)\nonumber\\
&\quad+\;\frac{2}{\hat{n}}s\,K_{\tfrac{\hat{n}-1}{2}}(C)\,\Phi_{1}\!\Big(\tfrac {\hat{n}}{2},1-\tfrac {\hat{n}}{2};\tfrac {\hat{n}}{2}+1;\ \kappa,\ Y\Big)
\Bigg],
\end{align}
where $K_{\nu}$ is the modified Bessel function of the second kind and
$\Phi_{1}$ is Humbert’s confluent hypergeometric function. Use $C>0$
in $K_{\nu}$ and in the power; the signed $sC$ appears only in $e^{sC}$
and in $Y$. By continuity, $F(0\mid\rho)=\tfrac12-\frac{\arcsin{\rho}}{\pi}$. Define
\[
B:=1-\rho^{2},\qquad C:=\frac{|z|}{B},\qquad s:=\sgn(z),\qquad
\kappa:=\frac{1-\rho}{2},\qquad Y:=-\,s\,C\,(1-\rho).
\]
\end{theorem}

The CDF, integral of the density first derived in \cite{NadarajahPogany2016} is,
\begin{equation}
\label{eq:pdf_prod}
F_Z(z\mid \rho) \;=\; \int_{-\infty}^{z}\frac{1}{\pi \sqrt{1-\rho^2}}
\exp\!\left(\frac{\rho Z}{1-\rho^2}\right)
K_0\!\left(\frac{|Z|}{1-\rho^2}\right)\, dZ.
\end{equation}
The CDF, integral of the density of the mean first derived in \cite{NadarajahPogany2016} is,
\begin{equation}
F_{\hat{Z}}(z) =
\int_{-\infty}^{z}\frac{ \hat{n}^{(\hat{n}+1)/2} 2^{(1-\hat{n})/2} \, |z|^{(\hat{n}-1)/2} }
{ \sqrt{\pi \, (1-\rho^2)} \, \Gamma(\hat{n}/2) }
\exp\!\left( \frac{\beta - \gamma}{2} z \right)
K_{\tfrac{1 - \hat{n}}{2}}\!\left( \frac{\beta + \gamma}{2} |z| \right)
\label{eq:pdf_prodMean}
\end{equation}
Using \eqref{eq:pdf_prodMean}, the integral of the density of the sum written,
\begin{equation}
F_{Z_{\Sigma}}(z) =
\int_{-\infty}^{z}\frac{ 2^{(1-\hat{n})/2} \, |z|^{(\hat{n}-1)/2} }
{ \sqrt{\pi \, (1-\rho^2)} \, \Gamma(\hat{n}/2) }
\exp\!\left( \frac{\hat{\beta} - \hat{\gamma}}{2} z \right)
K_{\tfrac{1 - \hat{n}}{2}}\!\left( \frac{\hat{\beta} + \hat{\gamma}}{2} |z| \right)
\label{eq:pdf_prodSum}
\end{equation}

\noindent
for $-\infty < z < \infty$, where
$\beta = \hat{n}/(1-\rho)$, $\gamma = \hat{n}/(1+\rho)$ and $\hat{\beta} = 1/(1-\rho)$, $\hat{\gamma} = 1/(1+\rho)$ 

Then, rearranging and expanding \eqref{eq:pdf_prod}, using $\gamma=\tfrac12$, it is possible to write the CDF
$F(z\mid\rho)=\Pr[Z\le z]$ of the product of correlated standard
Gaussians in terms of a normal--gamma mixture.
After the rescaling $y=(1-\rho^{2})\,g$,
\begin{align}
\label{eq:gammamix}
F(z\mid\rho)=\frac{1}{\sqrt{\pi}}\int_{0}^{\infty}
y^{-1/2}e^{-y}\,
\mathcal{N}\!\left(\frac{z-2\rho\,y}{\sqrt{2(1-\rho^{2})\,y}}\right)\,dy,    
\end{align}
For the case of sample mean, $\hat{Z}$ in \cite{NadarajahPogany2016} Theorem 2.2, which is \eqref{eq:pdf_prodMean} we will have the corresponding Normal-gamma mixture CDF,
\begin{align}
\label{eq:gammamixmean}
F(z\mid\rho)=\frac{1}{ \Gamma(\hat{n}/2)} \int_{0}^{\infty}
y^{\hat{n}/2-1}e^{-y}  \,
\mathcal{N}\!\left(\frac{\hat{n}z-2\rho\,y}{\sqrt{2(1-\rho^{2})\,y}}\right)\,dy,    
\end{align}
For the case of sample sum, $Z_{\Sigma}$  we will have the corresponding Normal-gamma mixture CDF,
\begin{align}
\label{eq:gammamixsum}
F(z\mid\rho)=\frac{1}{ \Gamma(\hat{n}/2)} \int_{0}^{\infty}
y^{\hat{n}/2-1}e^{-y}  \,
\mathcal{N}\!\left(\frac{z-2\rho\,y}{\sqrt{2(1-\rho^{2})\,y}}\right)\,dy,    
\end{align}

where $\mathcal{N}(\cdot)$ is the standard normal CDF. The integral \eqref{eq:gammamix}
equals the closed-form \eqref{eq:analytical_cdf_humbert} in Theorem \ref{thm:theorem1}, the integral \eqref{eq:gammamixmean} equals that of in Theorem \ref{thm:theorem2}, the integral \eqref{eq:gammamixsum} equals that of  in Theorem \ref{thm:theorem3} all pointwise.

\section{Proofs}
\subsection{Proof of Theorem 1}{\label{prf:thm1}}

Let's define the transformations following  \cite{MadanCarrChang1998}, \begin{align*}
n(v)&=\frac{z}{\sqrt{1-\rho^{2}}}\frac{1}{\sqrt{2(1-\rho^{2})+2 v^{2}}}, ~m(v)=\frac{v\sqrt{2}}{\sqrt{1-\rho^{2}}}\\
n&=\frac{z}  {\sqrt{2(1-\rho^{2})}}, ~m=\frac{\rho\sqrt{2}}{\sqrt{1-\rho^{2}}}.    
\end{align*}
Let's write the integral,

\begin{align}
F(z|\rho)
&= \int_{0}^{\infty} \Bigg[
\underbrace{\int_{-\infty}^{\rho}\mathcal{N}_n\!\left(\frac{z}{\sqrt{2(1-\rho^{2})y}}
+ \frac{\rho\sqrt{2 y} }{\sqrt{1-\rho^{2}}}\right)
y^{-1}e^{-y}\frac{1}{\Gamma(1/2)}\,dv}_{I_1}
\label{eq:I1gauss}
\\
&\qquad+
\underbrace{\int_{-\infty}^{\rho}\mathcal{N}_m\!\left(\frac{z}{\sqrt{2(1-\rho^{2})y}}
+ \frac{\rho\sqrt{2 y} }{\sqrt{1-\rho^{2}}}\right)
e^{-y}y^{0}\frac{1}{\Gamma(1/2)}\,dv}_{I_2}
\Bigg]\,dy\label{eq:I2gauss}
\end{align}
Then, we rewrite the equations \eqref{eq:I1gauss} and \eqref{eq:I2gauss}\ in terms of \cite{aslamzubeyir} form leading to Bessel function representations,
\begin{align*}
\phi\!\left(\frac{z-2\rho y}{\sqrt{2(1-\rho^2)y}}\right)
&= \frac{1}{\sqrt{2\pi}}
\exp\!\left(-\frac{(z-2\rho y)^2}{4(1-\rho^2)y}\right) \\
&= \frac{1}{\sqrt{2\pi}}
e^{\frac{\rho z}{1-\rho^2}}
\exp\!\left(-\frac{y}{1-\rho^2}-\frac{z^2}{4(1-\rho^2)y}\right),
\end{align*}
so that each $y$-integral becomes a standard Laplace--type integral
\[
\int_0^\infty t^{\nu-1}\exp\!\left(-\alpha t - \frac{\beta}{t}\right)dt
=2\left(\frac{\beta}{\alpha}\right)^{\nu/2} K_\nu\!\left(2\sqrt{\alpha\beta}\right),
\qquad \alpha,\beta>0.
\]

With $\alpha=\tfrac{1}{1-\rho^2}$ and $\beta=\tfrac{z^2}{4(1-\rho^2)}$, we have
\[
\frac{\beta}{\alpha}=\frac{z^2}{4},
\qquad
2\sqrt{\alpha\beta}=\frac{|z|}{1-\rho^2}.
\]

Hence, the two integrals become
\begin{align*}
I_1 &= \int_0^\infty 
\phi\!\left(\frac{z-2\rho y}{\sqrt{2(1-\rho^2)y}}\right) y^{-1}e^{-y}\,dy 
= \frac{2}{\sqrt{2\pi}}\,e^{\frac{\rho z}{1-\rho^2}}
K_0\!\left(\frac{|z|}{1-\rho^2}\right), \\
I_2 &= \int_0^\infty 
\phi\!\left(\frac{z-2\rho y}{\sqrt{2(1-\rho^2)y}}\right) e^{-y}\,dy 
= \frac{|z|}{\sqrt{2\pi}}\,e^{\frac{\rho z}{1-\rho^2}}
K_1\!\left(\frac{|z|}{1-\rho^2}\right).
\end{align*}
Now, we will make use of $n(v)$ and $m(v)$ to keep the Bessel functions fixed under interval $\left[-\infty,\rho\right]$ which we define as the support of $v$.
\begin{align}
I_2(v)m_v\, dv &= e^{n(v)m(v)} \left( \frac{\beta(v)}{\alpha(v)} \right)^{1/2}
K_{1}\left(\frac{|z|}{1-\rho^{2}}\right)\frac{\sqrt{2}}{\sqrt{1-\rho^{2}}}, \\
I_1(v)n_v\, dv &= -e^{n(v)m(v)} \left( \frac{\beta(v)}{\alpha(v)} \right)^{0}
K_{0}\left(\frac{|z|}{1-\rho^{2}}\right)\,
v\left(v^{2}+(1-\rho^{2}) \right)^{-3/2}
\end{align}
In terms of Bessel function elements we have,
\begin{align*}
\beta(v) &= \frac{z^{2}}{2(1-\rho^{2})}\left(\frac{1}{\sqrt{(1-\rho^{2})+v^{2}}}\right)^{2}, 
\qquad 
\alpha(v)=\frac{2 v^{2}}{1-\rho^{2}},\\
\frac{\beta(v)}{\alpha(v)}&=\left( \frac{z^{2}}{4\left(v^{2}+1-\rho^{2} \right)^{2}}\right)^{1/2},\\
u &= \frac{v}{v^{2}+(1-\rho^{2})}, ~v = \frac{\sqrt{1-\rho^{2}}\,u}{\sqrt{1-u^{2}}},~dv = \sqrt{1-\rho^{2}}\,(1-u^{2})^{-3/2}
\end{align*}
Then in terms of $u$,
\begin{align*}
\frac{\beta(u)}{\alpha(u)}&=\frac{1-u^{2}}{1-\rho^{2}},~
n_v\, dv = (1-\rho^{2})^{-1/2}u (1-u^{2})^{-1/2}, \\ 
~m_v\, dv &= \sqrt{2}(1-u^{2})^{-3/2},
~e^{m(v)n(v)} = e^{\frac{z\,u}{1-\rho^{2}}}.
\end{align*}
\begin{align}
\label{eq:I1}
I_2(u)m_u\, du &=2|z| 
e^{\frac{z u}{1-\rho^{2}}} 
\left( \frac{(1-u^{2})^{-1/2}}{2\left(1-\rho^{2}\right)}\right)
K_{1}\!\left(\frac{|z|}{1-\rho^{2}}\right)
\sqrt{2}, \\
\label{eq:I2}
I_1(u)n_u\, du &= 
-2\,e^{\frac{z u}{1-\rho^{2}}} 
K_{0}\!\left(\frac{|z|}{1-\rho^{2}}\right)
\,\frac{u}{\sqrt{1-u^{2}}} \frac{|z|}{\sqrt{2}(1-\rho^{2})}.
\end{align}

Then using \eqref{eq:I1} and \eqref{eq:I2} we can write
\begin{align*}
  I_1&=CK_{0}(C)\int_{-1}^{\rho}e^{\frac{zu}{1-\rho^{2}}}u (1-u^{2})^{-1/2}\frac{1}{\sqrt{2}} du\\
&=CK_{0}(C)e^{-C}\int_{0}^{1}e^{Cu(1+\rho)}\left(u(1+\rho)-1\right) (1-\left(u(1+\rho)-1)^{2}\right)^{-1/2}\sqrt{2} du\\
&=CK_{0}(C)e^{-C}\Biggl[\int_{0}^{1}e^{Cu(1+\rho)}\left(u^{-1/2}(1+\rho)^{1/2}\right) \left(1-\left(\frac{u(1+\rho)}{2}\right)^{-1/2}\right) du\\
&-\int_{0}^{1}e^{Cu(1+\rho)}\left(u^{1/2}(1+\rho)^{3/2}\right) \left(1-\left(\frac{u(1+\rho)}{2}\right)\right)^{-1/2} du\Biggr]
\end{align*}

\begin{align*}
  I_2&=CK_{1}(C)\int_{-1}^{\rho}e^{\frac{zu}{1-\rho^{2}}} (1-u^{2})^{-1/2}\sqrt{2}du\\
&=CK_{1}(C)e^{-C}\int_{0}^{1}e^{Cu(1+\rho)}(1-\left(u(1+\rho)-1)^{2}\right)^{-1/2}\sqrt{2} du(1+\rho)\\
&=CK_{1}(C)e^{-C}\int_{0}^{1}e^{Cu(1+\rho)}u^{-1/2}(1+\rho)^{1/2} \left(1-\left(\frac{u(1+\rho)}{2}\right)\right)^{-1/2} du
\end{align*}

Then using Humbert function  \cite{Humbert1920} or  \cite{Gradshteyn2014}, 
\begin{equation}
\Phi_{1}(\alpha,\beta;\gamma;x,y)
= \frac{\Gamma(\gamma)}{\Gamma(\alpha)\Gamma(\gamma-\alpha)}
\int_{0}^{1} u^{\alpha-1}(1-u)^{\gamma-\alpha-1}
(1-yu)^{-\beta}e^{xu}\,du.
\label{eq:HumbertIntegral}
\end{equation}
\begin{align}
\int_{0}^{1} e^{C u(1+\rho)}\,u^{1/2}
\left(1 - \frac{u(1+\rho)}{2}\right)^{-1/2} \, du
&=
\frac{\Phi_{1}\!\left(\tfrac{1}{2},\,\tfrac{1}{2};\,\tfrac{3}{2};\,\tfrac{1+\rho}{2},\,\operatorname{sgn}(z)\,C(1+\rho)\right)}{\frac{1}{2}},
\\[4pt]
\int_{0}^{1} e^{C u(1+\rho)}\,u^{1/2}
\left(1 - \frac{u(1+\rho)}{2}\right)^{-1/2} \, du
&= \frac{
\Phi_{1}\!\left(\tfrac{3}{2},\,\tfrac{1}{2};\,\tfrac{5}{2};\,\tfrac{1+\rho}{2},\,\operatorname{sgn}(z)\,C(1+\rho)\right)}{{\frac{3}{2}}}.
\label{eq:Phi1Integral}
\end{align}
The constants outside Humbert functions come from the fact that $\Gamma(1+\alpha)=\Gamma(\alpha)\alpha$ and since there is no $(1+u)$ term we set $\gamma=\alpha+1$. Moreover, due to these two conditions $\Gamma(x)$ terms cancel out and therefore, we are left with only $\alpha$ term in all Humbert form integrals. 

Finally, collecting all constants and special functions yields  equation \eqref{eq:analytical_cdf_humbert} in Theorem \ref{thm:theorem1}.
\subsection{Proof of Theorem 2}{\label{prf:thm2}}
Let's define the transformations similar to Proof  \ref{prf:thm1} considering the random variable defined in Theorem \ref{thm:theorem2},
\begin{align*}
n(v)&=\frac{\hat{n}z}{\sqrt{1-\rho^{2}}}\frac{1}{\sqrt{2(1-\rho^{2})+2 v^{2}}}, ~m(v)=\frac{v\sqrt{2}}{\sqrt{1-\rho^{2}}}\\
n&=\frac{\hat{n} z}  {\sqrt{2(1-\rho^{2})}}, ~m=\frac{\rho\sqrt{2}}{\sqrt{1-\rho^{2}}}.    
\end{align*}

Let's write the integral,

\begin{align}
F(z|\rho)
&= \int_{0}^{\infty} \Bigg[
\underbrace{\int_{-\infty}^{\rho}\mathcal{N}_n\!\left(\frac{\hat{n} z}{\sqrt{2(1-\rho^{2})y}}
+ \frac{\rho\sqrt{2 y} }{\sqrt{1-\rho^{2}}}\right)
y^{\hat{n}-1}e^{-y}\frac{1}{\Gamma(\hat{n}/2)}\,dv}_{I_1}
\label{eq:I1gaussmean}
\\
&\qquad+
\underbrace{\int_{-\infty}^{\rho}\mathcal{N}_m\!\left(\frac{z}{\sqrt{2(1-\rho^{2})y}}
+ \frac{\rho\sqrt{2 y} }{\sqrt{1-\rho^{2}}}\right)
e^{-y}y^{\hat{n}}\frac{1}{\Gamma(\hat{n}/2)}\,dv}_{I_2}
\Bigg]\,dy\label{eq:I2gaussmean}
\end{align}
Then, we rewrite the equations \eqref{eq:I1gauss} and \eqref{eq:I2gauss}\ in terms of \cite{aslamzubeyir} form leading to Bessel function representations,
\begin{align*}
\phi\!\left(\frac{\hat{n}z-2\rho y}{\sqrt{2(1-\rho^2)y}}\right)
&= \frac{1}{\sqrt{2\pi}}
\exp\!\left(-\frac{(\hat{n}z-2\rho y)^2}{4(1-\rho^2)y}\right) \\
&= \frac{1}{\sqrt{2\pi}}
e^{\frac{\rho z}{1-\rho^2}}
\exp\!\left(-\frac{y}{1-\rho^2}-\frac{\hat{n}^{2}z^2}{4(1-\rho^2)y}\right),
\end{align*}
so that each $y$-integral becomes a standard Laplace--type integral
\[
\int_0^\infty t^{\nu-1}\exp\!\left(-\alpha t - \frac{\beta}{t}\right)dt
=2\left(\frac{\beta}{\alpha}\right)^{\nu/2} K_\nu\!\left(2\sqrt{\alpha\beta}\right),
\qquad \alpha,\beta>0.
\]

With $\alpha=\tfrac{1}{1-\rho^2}$ and $\beta=\tfrac{z^2}{4(1-\rho^2)}$, we have
\[
\frac{\beta}{\alpha}=\frac{z^2}{4},
\qquad
2\sqrt{\alpha\beta}=\frac{|z|}{1-\rho^2}.
\]

Hence, the two integrals become
\begin{align*}
I_1 &= \int_0^\infty 
\phi\!\left(\frac{\hat{n}z-2\rho y}{\sqrt{2(1-\rho^2)y}}\right) y^{\tfrac{\hat{n}-3}{2}}e^{-y}\,dy 
=  2 \left(\frac{\beta}{\alpha}\right)^{\tfrac{\hat{n}-1}{4}}{\sqrt{2\pi}}\,e^{\frac{\rho z}{1-\rho^2}}
K_{\tfrac{\hat{n}-1}{2}}\!\left(\frac{|z|}{1-\rho^2}\right), \\
I_2 &= \int_0^\infty 
\phi\!\left(\frac{z-2\rho y}{\sqrt{2(1-\rho^2)y}}\right) y^{\hat{n}-1}e^{-y}\,dy 
= 2 \left(\frac{\beta}{\alpha}\right)^{\tfrac{\hat{n}+1}{4}} {\sqrt{2\pi}}\,e^{\frac{\rho z}{1-\rho^2}}
K_{\tfrac{\hat{n}+1}{2}}\!\left(\frac{|z|}{1-\rho^2}\right), \\
\end{align*}
Now, we will make use of $n(v)$ and $m(v)$ to keep the Bessel functions fixed under interval $\left[-\infty,\rho\right]$ which we define as the support of $v$.
\begin{align}
I_2(v)m_v\, dv &= e^{n(v)m(v)} \left( \frac{\beta(v)}{\alpha(v)} \right)^{\tfrac{\hat{n}+1}{4}
}
K_{\frac{\hat{n}+1}{2}}\left(\frac{|z|}{1-\rho^{2}}\right)\frac{\sqrt{2}}{\sqrt{1-\rho^{2}}}, \\
I_1(v)n_v\, dv &= -e^{n(v)m(v)} \left( \frac{\beta(v)}{\alpha(v)} \right)^{\tfrac{\hat{n}-1}{4}}
K_{\tfrac{\hat{n}-1}{2}}\left(\frac{|z|}{1-\rho^{2}}\right)\,
v\left(v^{2}+(1-\rho^{2}) \right)^{-3/2}
\end{align}
For the rest we proceed with the same simplifications and boundary adjustments and variable transformations ($v\rightarrow{u}$), finally evaluation of the integral, obtained after these operations, in terms of Humbert's confluent function definition in proof \ref{prf:thm1} to match necessary parameters, then the result in Theorem \ref{thm:theorem2} follows.

\subsection{Proof of Theorem 3}
Using the transformations below,
\begin{align*}
n(v)&=\frac{z}{\sqrt{1-\rho^{2}}}\frac{1}{\sqrt{2(1-\rho^{2})+2 v^{2}}}, ~m(v)=\frac{v\sqrt{2}}{\sqrt{1-\rho^{2}}}\\
n&=\frac{z}  {\sqrt{2(1-\rho^{2})}}, ~m=\frac{\rho\sqrt{2}}{\sqrt{1-\rho^{2}}}.    
\end{align*} 
then considering the Normal-gamma mixture representation for sum $Z_{\Sigma}$ in equation \eqref{eq:pdf_prodSum},
then following Proofs \ref{prf:thm1} and \ref{prf:thm2} for the rest the result in Theorem \ref{thm:theorem3} follows.


\section{Figures}\label{sec6}

\begin{figure}[h]
\centering
\includegraphics[width=0.9\textwidth]{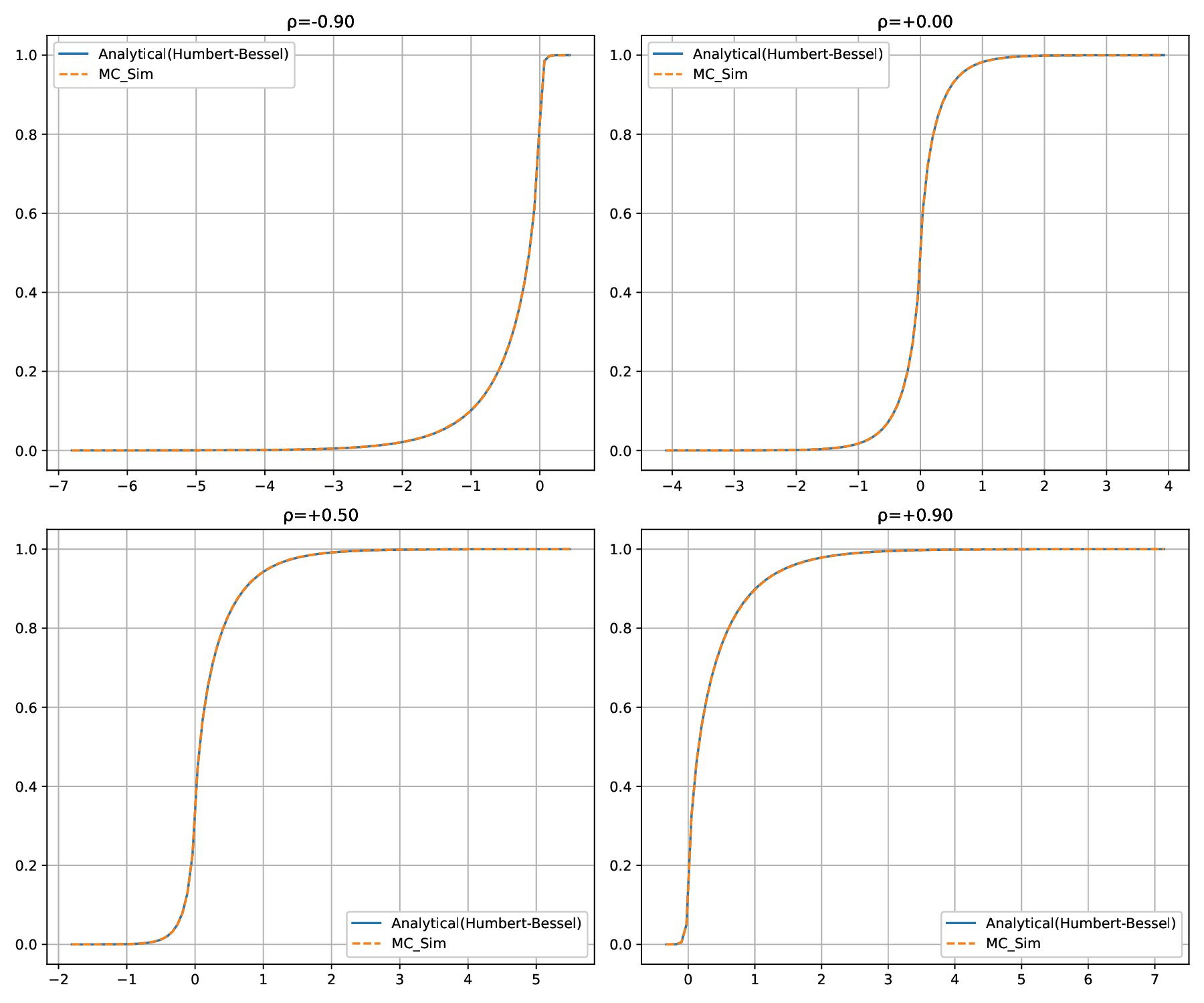}
\caption{Analytical CDF vs.\ MC empirical CDF (Humbert-Bessel vs MC).}
\label{fig:cdf_50_rhoMC}
\end{figure}

\begin{figure}[h]
\centering
\includegraphics[width=0.9\textwidth]{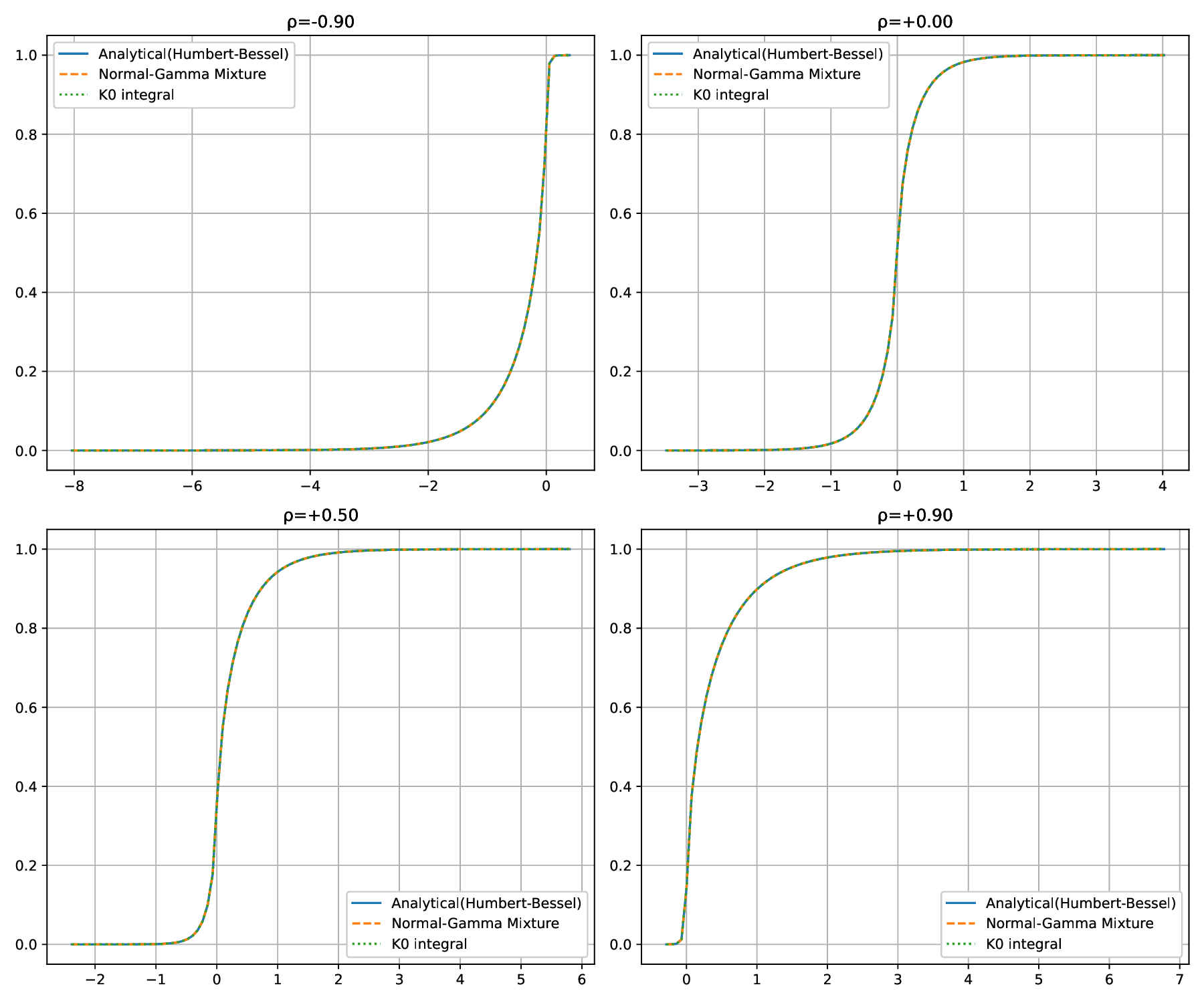}
\caption{CDF in three methods at $\rho=-0.9,0.0,0.5,0.9$.}
\label{fig:cdf_50_rho}
\end{figure}

\begin{figure}[h]
\centering
\includegraphics[width=0.9\textwidth]{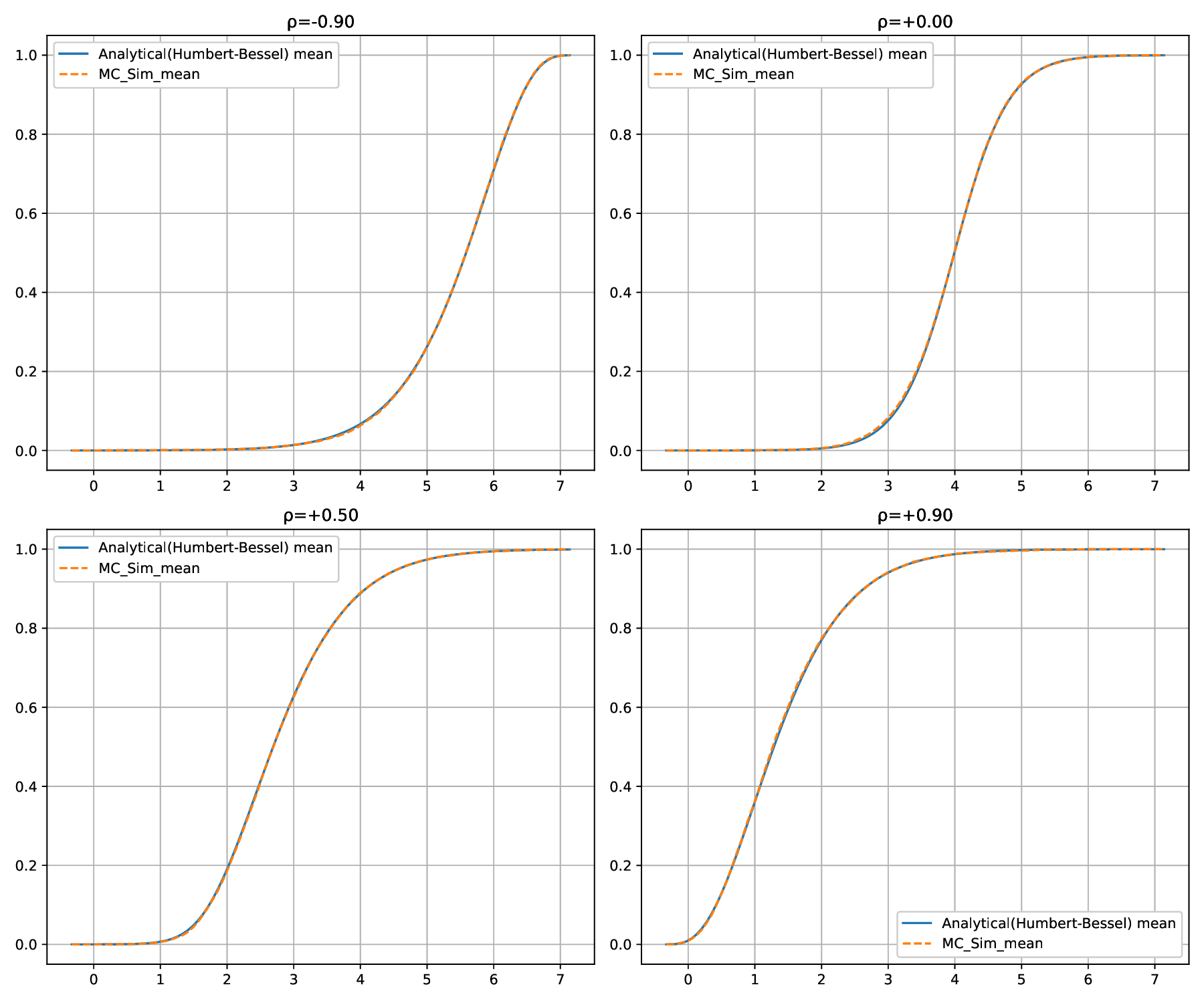}
\caption{Analytical CDF vs.\ MC empirical mean CDF (Humbert-Bessel vs MC).}
\label{fig:cdf_50_rhoMCMean}
\end{figure}

\begin{figure}[h]
\centering
\includegraphics[width=0.9\textwidth]{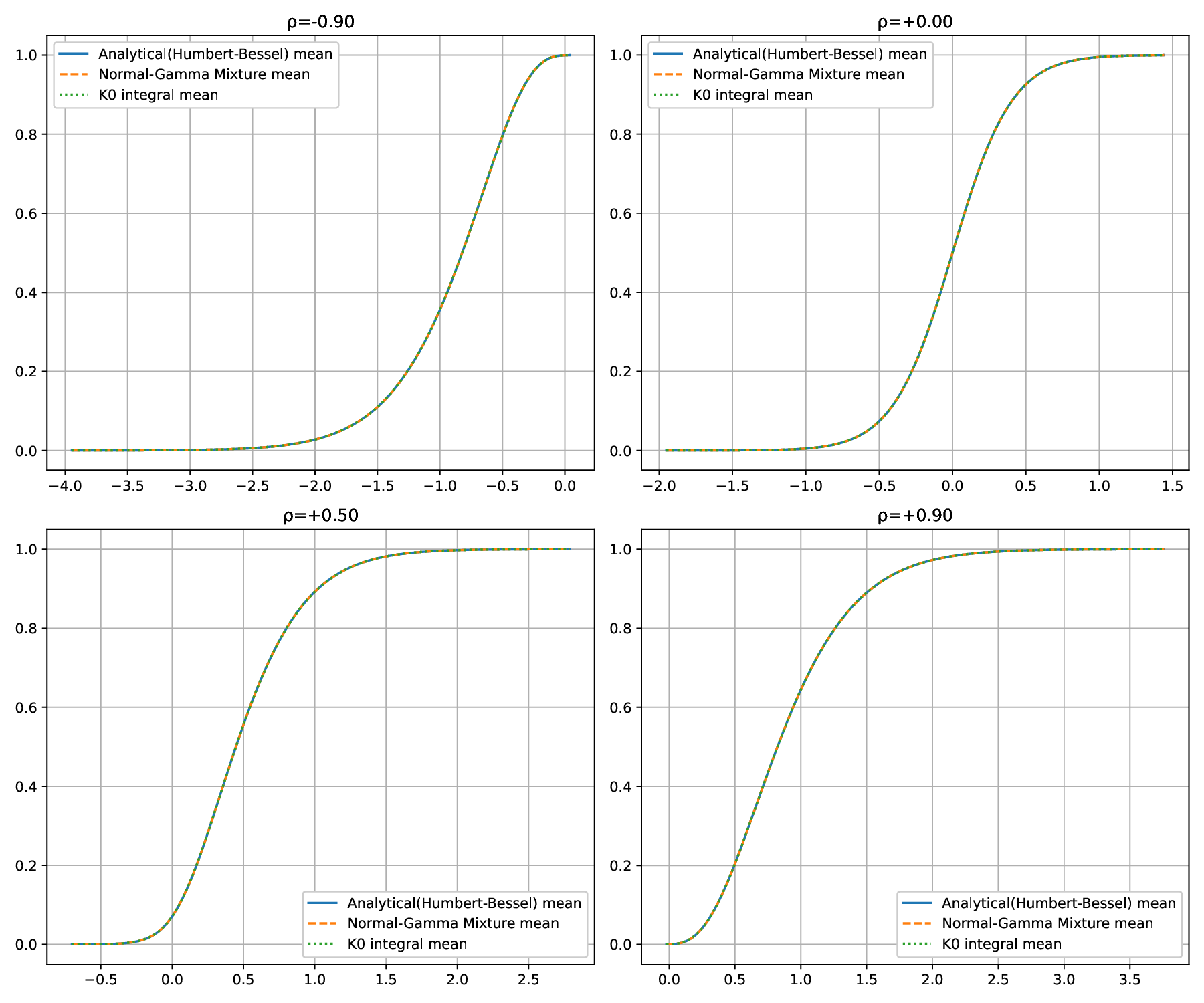}
\caption{Mean CDF in three methods at $\rho=-0.9,0.0,0.5,0.9$.}
\label{fig:cdf_50_rhoMean}
\end{figure}

\begin{figure}[h]
\centering
\includegraphics[width=0.9\textwidth]{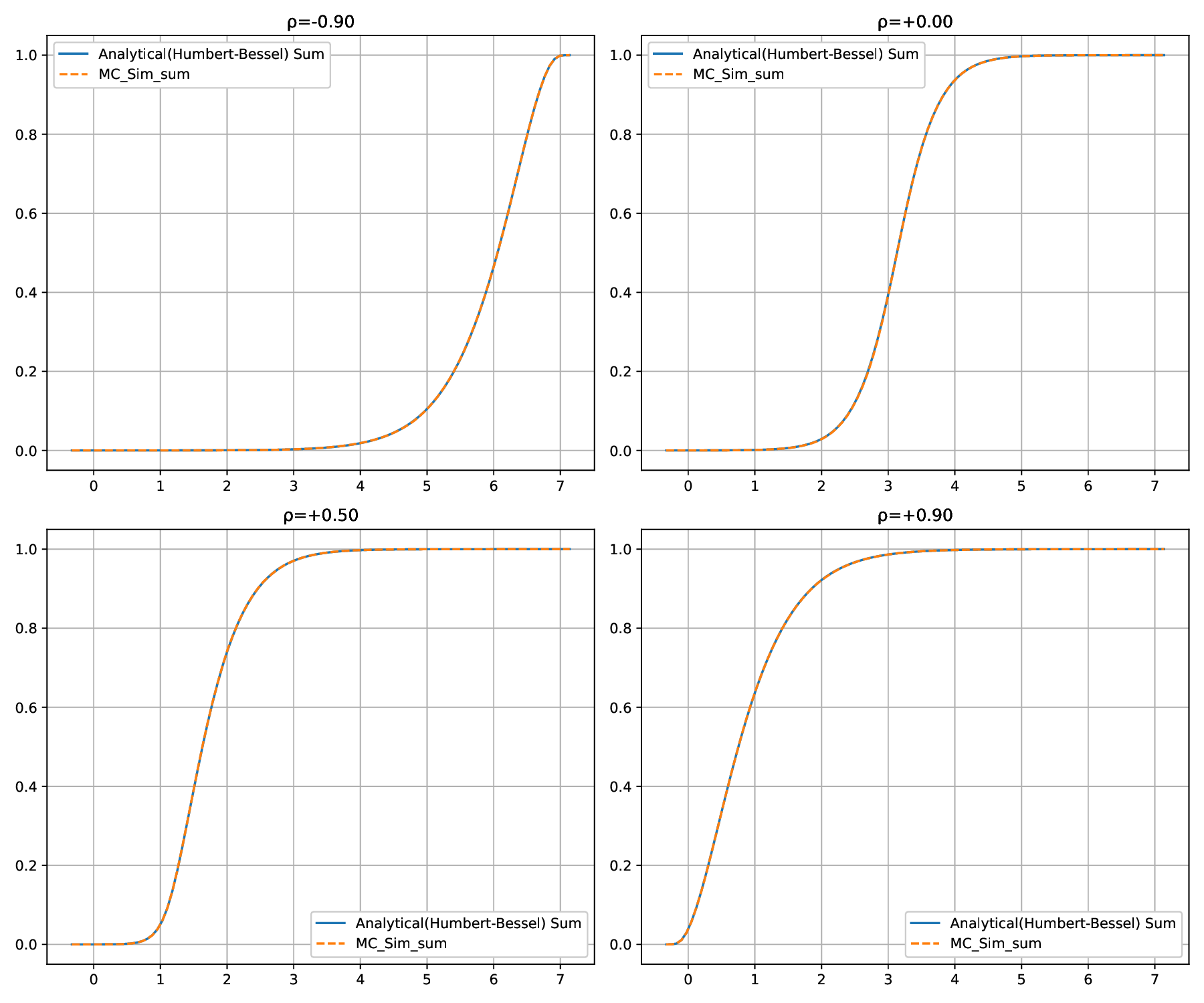}
\caption{Analytical CDF vs.\ MC empirical sum CDF (Humbert-Bessel vs MC).}
\label{fig:cdf_50_rhoMCSum}
\end{figure}

\begin{figure}[h]
\centering
\includegraphics[width=0.9\textwidth]{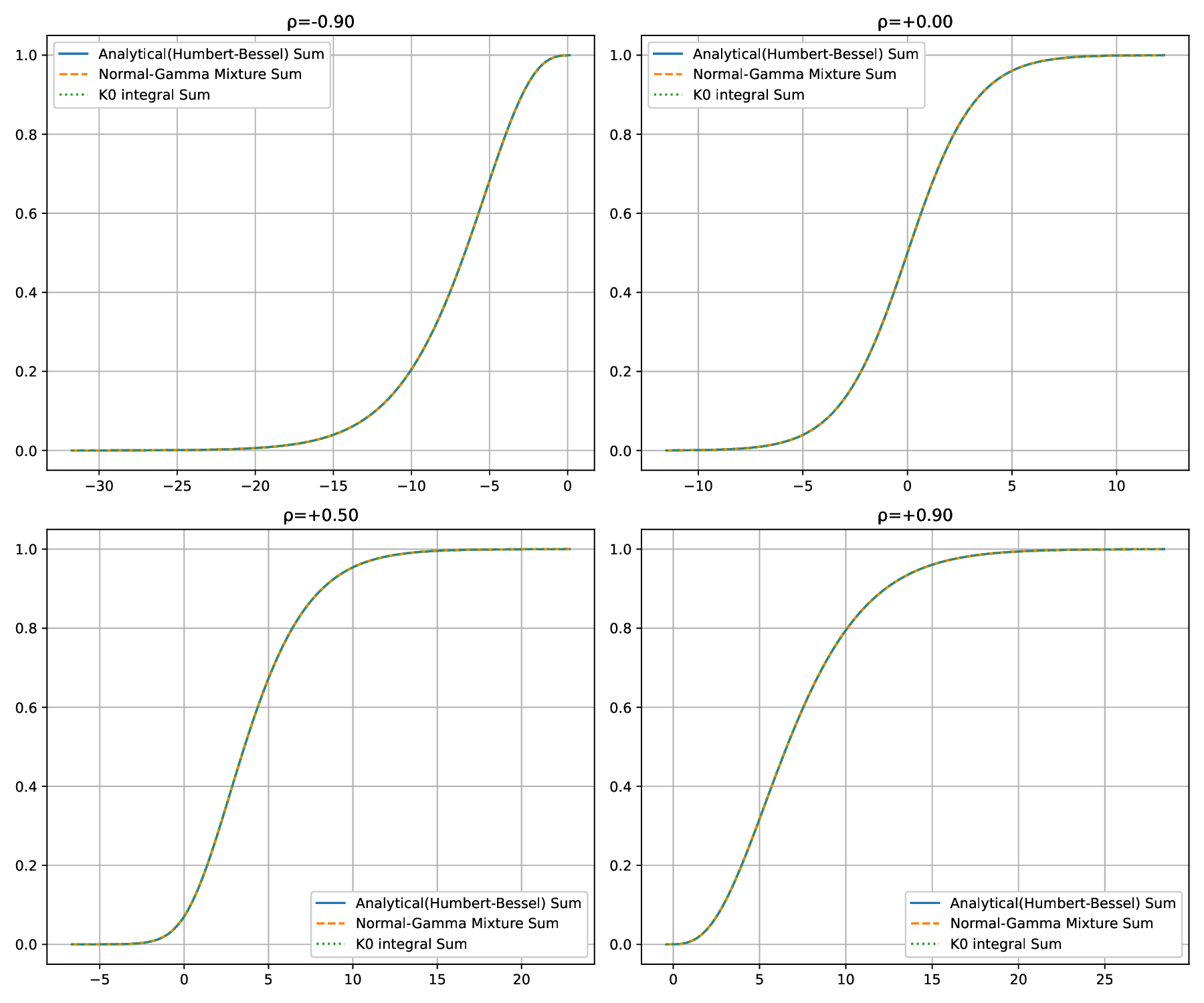}
\caption{Sum CDF in three methods at $\rho=-0.9,0.0,0.5,0.9$.}
\label{fig:cdf_50_rhoSum}
\end{figure}

The numerical experiments confirm the formula's validity at machine precision. 
Figure~\ref{fig:cdf_50_rho} shows that the numerical integral of the analytical PDF with 
$K_0$ (the modified Bessel function of the second kind at $\nu=0$) coincides exactly with 
the closed-form CDF in~\eqref{eq:analytical_cdf_humbert}. 
Moreover, Figure~\ref{fig:cdf_50_rhoMC} demonstrates perfect agreement with the 
Monte Carlo CDF of the normal-product variable, confirming the robustness of 
\eqref{eq:analytical_cdf_humbert}. In figures, \ref{fig:cdf_50_rhoMCMean} and \ref{fig:cdf_50_rhoMCSum} we see again perfect alignment  with  MC based CDF. In figures \ref{fig:cdf_50_rhoMean} and \ref{fig:cdf_50_rhoSum} we again see equality at machine precision level. Therefore, integral of densities in \eqref{eq:pdf_prodMean} and \eqref{eq:pdf_prodSum}. Moreover, Normal-Gamma mixture of the mean in equation \eqref{eq:gammamixmean} and Normal-Gamma mixture of the sum in equation \eqref{eq:gammamixsum} show exact alignment with \eqref{eq:analytical_cdf_humbertMean} and \eqref{eq:analytical_cdf_humbertSum} respectively at machine precision level.

\section{Tables}

Regarding performance, Table~\ref{tab:TablePerformance} shows that the closed-form 
formula in~\eqref{eq:analytical_cdf_humbert} achieves \textbf{2--3$\times$} higher 
efficiency than \eqref{eq:pdf_prod} (numerical integration of the density) and Normal-gamma mixture, \eqref{eq:gammamix} methods while maintaining the same machine-precision 
accuracy. 
In Table \ref{tab:TablePerformanceMean}, the distribution function of mean, we see similar stability plus an even better precision and computation performance (\textbf{3--10$\times$} more efficient). We observe similar features in Table \ref{tab:TablePerformanceSum} for the cumulative distribution of the sum as well. Therefore, we can confirm that the equation \eqref{eq:analytical_cdf_humbertMean} and \eqref{eq:analytical_cdf_humbertSum} work quite accurately to compute CDF of both product Gaussian mean and product Gaussian  sum random variables.

\begin{table}[h!]
\centering
\caption{Benchmark timings (seconds) and max/mean differences across methods for selected $\rho$.}
\label{tab:TablePerformance}
\begin{tabular}{lcccc}
\toprule
\textbf{Method} & \(\rho=-0.9\) & \(\rho=0.0\) & \(\rho=+0.5\) & \(\rho=+0.9\) \\
\midrule
Analytical (Humbert–Bessel) & 0.044 s & 0.044 s & 0.043 s & 0.044 s \\
Normal–Gamma Mixture        & 0.081 s & 0.075 s & 0.071 s & 0.059 s \\
K0 integral                 & 0.140 s & 0.128 s & 0.132 s & 0.134 s \\
\midrule
$\max|\text{Analytical}-\text{Mixture}|$ & $6.128\times 10^{-14}$ & $3.187\times 10^{-13}$ & $4.907\times 10^{-14}$ & $4.108\times 10^{-14}$ \\
Mean                                       & $7.845\times 10^{-15}$ & $1.155\times 10^{-14}$ & $5.136\times 10^{-15}$ & $7.241\times 10^{-15}$ \\
\midrule
$\max|\text{Analytical}-\text{K0}|$        & $1.241\times 10^{-07}$ & $5.410\times 10^{-08}$ & $6.247\times 10^{-08}$ & $1.241\times 10^{-07}$ \\
Mean                                       & $3.418\times 10^{-09}$ & $1.651\times 10^{-09}$ & $1.777\times 10^{-09}$ & $3.418\times 10^{-09}$ \\
\midrule
$\max|\text{Mixture}-\text{K0}|$           & $1.241\times 10^{-07}$ & $5.410\times 10^{-08}$ & $6.247\times 10^{-08}$ & $1.241\times 10^{-07}$ \\
Mean                                       & $3.418\times 10^{-09}$ & $1.651\times 10^{-09}$ & $1.777\times 10^{-09}$ & $3.418\times 10^{-09}$ \\
\bottomrule
\end{tabular}
\end{table}
\begin{table}[ht]
\centering
\caption{Benchmark timings (seconds) and max/mean differences across methods for selected $\rho$ (mean $\hat{Z}$ variable).}
\label{tab:TablePerformanceMean}
\begin{tabular}{lcccc}
\toprule
\textbf{Method} & $\rho=-0.9$ & $\rho=0.0$ & $\rho=+0.5$ & $\rho=+0.9$ \\
\midrule
Analytical (Humbert--Bessel) & 0.038 s & 0.009 s & 0.010 s & 0.011 s \\
Normal--Gamma Mixture        & 0.104 s & 0.079 s & 0.071 s & 0.087 s \\
K0 integral                  & 0.117 s & 0.050 s & 0.083 s & 0.103 s \\
\midrule
max $|\text{Analytical}-\text{Mixture}|$ & $6.839\times 10^{-14}$ & $7.794\times 10^{-14}$ & $1.085\times 10^{-13}$ & $4.172\times 10^{-15}$ \\
Mean                                       & $1.755\times 10^{-15}$ & $2.894\times 10^{-15}$ & $2.311\times 10^{-15}$ & $2.681\times 10^{-16}$ \\
\midrule
max $|\text{Analytical}-\text{K0}|$        & $2.875\times 10^{-11}$ & $4.258\times 10^{-11}$ & $1.343\times 10^{-10}$ & $6.711\times 10^{-10}$ \\
Mean                                       & $3.215\times 10^{-12}$ & $1.858\times 10^{-12}$ & $4.008\times 10^{-12}$ & $1.552\times 10^{-11}$ \\
\midrule
max $|\text{Mixture}-\text{K0}|$           & $2.875\times 10^{-11}$ & $4.258\times 10^{-11}$ & $1.343\times 10^{-10}$ & $6.711\times 10^{-10}$ \\
Mean                                       & $3.217\times 10^{-12}$ & $1.860\times 10^{-12}$ & $4.010\times 10^{-12}$ & $1.552\times 10^{-11}$ \\
\bottomrule
\end{tabular}
\end{table}

\begin{table}[ht]
\centering
\caption{Benchmark timings (seconds) and max/mean differences across methods for selected $\rho$ (sum $Z_{\Sigma}$ variable).}
\label{tab:TablePerformanceSum}
\begin{tabular}{lcccc}
\toprule
\textbf{Method} & $\rho=-0.9$ & $\rho=0.0$ & $\rho=+0.5$ & $\rho=+0.9$ \\
\midrule
Analytical (Humbert--Bessel) & 0.039 s & 0.010 s & 0.011 s & 0.012 s \\
Normal--Gamma Mixture        & 0.112 s & 0.084 s & 0.072 s & 0.094 s \\
K0 integral                  & 0.123 s & 0.076 s & 0.094 s & 0.136 s \\
\midrule
max $|\text{Analytical}-\text{Mixture}|$ & $2.021\times 10^{-14}$ & $7.122\times 10^{-14}$ & $4.257\times 10^{-14}$ & $8.793\times 10^{-15}$ \\
Mean                                       & $1.048\times 10^{-15}$ & $2.437\times 10^{-15}$ & $1.410\times 10^{-15}$ & $3.654\times 10^{-16}$ \\
\midrule
max $|\text{Analytical}-\text{K0}|$        & $1.350\times 10^{-09}$ & $4.488\times 10^{-11}$ & $1.099\times 10^{-09}$ & $8.669\times 10^{-11}$ \\
Mean                                       & $2.651\times 10^{-11}$ & $3.765\times 10^{-12}$ & $1.451\times 10^{-11}$ & $3.776\times 10^{-12}$ \\
\midrule
max $|\text{Mixture}-\text{K0}|$           & $1.350\times 10^{-09}$ & $4.488\times 10^{-11}$ & $1.099\times 10^{-09}$ & $8.669\times 10^{-11}$ \\
Mean                                       & $2.651\times 10^{-11}$ & $3.767\times 10^{-12}$ & $1.451\times 10^{-11}$ & $3.776\times 10^{-12}$ \\
\bottomrule
\end{tabular}
\end{table}

\section{Conclusion}

In this paper, we have derived a one-line closed-form expression for the cumulative distribution function (CDF) of the product of zero-mean correlated Gaussian random variables. The final formula, expressed in terms of Humbert's confluent hypergeometric function and the modified Bessel function of the second kind, provides a compact and analytically tractable representation.  

Comprehensive numerical experiments and Monte Carlo simulations confirm the accuracy and computational efficiency of the proposed formula across a wide range of parameter settings. Owing to its closed-form nature and high precision, the result offers a practical tool for applications in modeling non-linear signals, quantitative finance, probability theory, and related areas.  

Future research will focus on extending this methodology to the case of Gaussian variables with non-zero means and exploring potential generalizations to higher-dimensional settings.

\clearpage
\bibliography{sn-bibliography}
\end{document}